\documentclass[fleqn, leqno, 12pt]{article}

\usepackage{amsmath, amssymb, bm}
\usepackage[dvips]{graphicx, color}

\newcommand{\remcomm}{Remark 2.1}
\newcommand{\propmain}{Proposition 2.1}
\newcommand{\thmmain}{Theorem 2.1}

\newcommand{\thmbasic}{Theorem 3.1}

\newcommand{\propisom}{Proposition 3.1}
\newcommand{\propkernel}{Proposition 3.2}
\newcommand{\proppsd}{Proposition 3.3}
\newcommand{\lemwel}{Lemma 3.1}
\newcommand{\propvzc}{Proposition 3.4}
\newcommand{\propab}{Proposition 3.5}

\newcommand{\lemsi}{Lemma 4.1}
\newcommand{\lemti}{Lemma 4.2}

\newcommand{\propqua}{Proposition 5.1}
\newcommand{\remqua}{Remark 5.1}
\newcommand{\defper}{Definition 5.1}
\newcommand{\lemhone}{Lemma 5.1}
\newcommand{\lemgood}{Lemma 5.2}
\newcommand{\lemcont}{Lemma 5.3}
\newcommand{\lemqua}{Lemma 5.4}
\newcommand{\lemfs}{Lemma 5.5}
\newcommand{\lemaa}{Lemma 5.6}

\newcommand{\Bone}{B1}
\newcommand{\Btwo}{B2}
\newcommand{\Bthree}{B3}

\newcommand{\Kone}{K1}

\newcommand{\eqpde}{1.1}

\newcommand{\eqholder}{2.1}
\newcommand{\eqtpe}{2.2}
\newcommand{\Hone}{H1}
\newcommand{\HoneFirst}{H1-1}
\newcommand{\HoneSecond}{H1-2}
\newcommand{\HoneThird}{H1-3}
\newcommand{\HoneFourth}{H1-4}
\newcommand{\Htwo}{H2}
\newcommand{\Hthree}{H3}
\newcommand{\Hfour}{H4}
\newcommand{\Hfive}{H5}

\newcommand{\eqH}{3.1}
\newcommand{\eqThe}{3.2}
\newcommand{\eqgi}{3.3}
\newcommand{\eqlon}{3.4}
\newcommand{\eqtoa}{3.5}
\newcommand{\eqria}{3.6}
\newcommand{\eqrsc}{3.7}
\newcommand{\eqpco}{3.8}
\newcommand{\eqvarp}{3.9}
\newcommand{\eqPper}{3.10}
\newcommand{\eqEq}{3.11}

\newcommand{\eqfrac}{4.1}
\newcommand{\eqxn}{4.2}
\newcommand{\eqyn}{4.3}
\newcommand{\eqmapdh}{4.4}
\newcommand{\eqgam}{4.5}
\newcommand{\eqdefa}{4.6}
\newcommand{\eqmuz}{4.7}
\newcommand{\eqfly}{4.8}
\newcommand{\eqsvec}{4.9}
\newcommand{\eqrto}{4.10}
\newcommand{\eqluo}{4.11}
\newcommand{\eqpin}{4.12}
\newcommand{\eqwg}{4.13}
\newcommand{\eqtuv}{4.14}
\newcommand{\eqwta}{4.15}
\newcommand{\eqzta}{4.16}
\newcommand{\eqcd}{4.17}
\newcommand{\eqcdt}{4.18}

\newcommand{\eqhc}{5.1}
\newcommand{\eqhu}{5.2}
\newcommand{\eqho}{5.3}
\newcommand{\eqau}{5.4}
\newcommand{\eqax}{5.5}
\newcommand{\eqg}{5.6}
\newcommand{\eqdva}{5.7}
\newcommand{\eqdtv}{5.8}
\newcommand{\eqdt}{5.9}
\newcommand{\eqdh}{5.10}

\newcommand{\refAma}{A}
\newcommand{\refBKST}{BKST}
\newcommand{\refABB}{ABB}
\newcommand{\refCRbif}{CR1}
\newcommand{\refCRhopf}{CR}

\newcommand{\refSpectral}{EE}
\newcommand{\refGMW}{GMW}
\newcommand{\refKsym}{K1}
\newcommand{\refKcom}{K2}
\newcommand{\refKsome}{K3}
\newcommand{\refKhi}{K4}
\newcommand{\refKiel}{Ki}
\newcommand{\refliu}{LMR}
\newcommand{\refLiZY}{LiZY}
\newcommand{\refMS}{MS}

\newcommand{\refWYZ}{WYZ}

\newcommand{\ubold}{{\bm u}}
\newcommand{\hbold}{{\bm h}}

\newcommand{\ds}{\displaystyle}

\newcommand{\Uc}{U_{\hbox{c}}}
\newcommand{\Vc}{V_{\hbox{c}}}
\newcommand{\Ac}{A_{\hbox{c}}}
\newcommand{\us}{u_\star}

\newcommand{\eone}{{\bm e}_1}

\newcommand{\psis}{\psi_\star}
\newcommand{\psish}{\psi_\sharp}

\newcommand{\Vst}{V_\star}
\newcommand{\Vsh}{V_\sharp}
\newcommand{\Ast}{A_\star}
\newcommand{\Ash}{A_\sharp}

\newcommand{\cA}{ {\cal A}}
\newcommand{\cL}{ {\cal L}}
\newcommand{\cU}{ {\cal U}}
\newcommand{\cV}{ {\cal V} }

\newcommand{\cX}{ {\cal X} }
\newcommand{\cY}{ {\cal Y} }
\newcommand{\cZ}{ {\cal Z} }

\newcommand{\spnc}{\spn_{\hbox{c}}}

\newcommand{\sech}{\mathrm{sech}\,}

\newcommand{\re}{\mathrm{Re\,}}
\newcommand{\im}{\mathrm{Im\,}}

\newcommand{\linear}{ {\cal L} }
\newcommand{\domain}{ {\cal D} }
\newcommand{\range}{ {\cal R} }
\newcommand{\kernel}{ {\cal N} }

\newcommand{\spn}{\mathrm{span}}

\newcommand{\directsum}{\oplus}

\newcommand{\codim}{\mathrm{codim}\,}

\newcommand{\n}{\par\noindent}

\newcommand{\propspace}{\vspace{2.1ex}\n} 

\newcommand{\eqspace}{\vspace{1.3ex}\n}
\newcommand{\newchapter}{\vspace{7.0 ex}}  

\newcommand{\q}{\quad}

\newcommand{\medn}{\medskip\n}

\newcommand{\npage}{\vfil \break}
 \newcommand{\End}{\hfill $\Box$}

\newcommand{\moji}[1] {\hspace{0.7 em}\text{#1}\hspace {0.7 em}} 

\newcommand{\semicolon}{\,;\,}
\newcommand{\mycolon}{\,:\,}

\newcommand{\N}{ \mathbb{N} }
\newcommand{\R}{ \mathbb{R} }
\newcommand{\Z}{ \mathbb{Z} }
\newcommand{\C}{\mathbb{C}}

\newcommand{\e}{\varepsilon}

\newcommand{\lam}{\lambda}
\newcommand{\Lam}{\Lambda}

\begin{document}

\setlength{\baselineskip}{6.5 mm}


\centerline{ \bf  The Hopf bifurcation theorem in Banach spaces}

  \vskip 0.3 cm

 \centerline{Tadashi KAWANAGO}
\vskip 0.3 cm
\centerline{Faculty of Education, Saga University, Saga 840-8502, Japan}
\centerline{e-mail: tadashi@cc.saga-u.ac.jp}

\vspace{1.5 cm}

{\bf Abstract}
\n                                                             
We prove a Hopf bifurcation theorem in general Banach spaces, which   
improves a  classical result by Crandall and Rabinowitz. Actually, 
our theorem does not need any compactness conditions, which leads to 
wider applications. In particular, our  theorem can be applied to
semilinear and quasi-linear partial differential equations 
in unbounded domains of $\R^n$.

\newchapter

\centerline{\bf  1. Introduction}
Concerning the Hopf bifurcation theorems in infinite dimensions,
a lot of versions have been proved until now
 (see e.g.   [\refCRhopf],  [\refAma],  [\refliu],
[\refGMW] and the references therein).
Among them [\refCRhopf, Theorem 1.11] by Crandall and Rabinowitz is 
one of most important results. It is a theorem for
abstract  semilinear equations
and has been well applied so far
to various studies
because of its generality (see e.g. [\refGMW] and [\refWYZ]). 
 It needs, however,
  some compactness condition, and, consequently, can not
  be applied   to partial differential equations
  in unbouded domains  of $\R^n$.
 
 On the other hand, Hopf bifurcation in partial differential equations 
  in the unbounded domain of $\R^n$ has been  studied 
  more recently and
   Hopf bifurcation theorems applicable to
 such studies were proven
(see e.g. [\refLiZY], [\refMS] and [\refBKST]). 
As far as the author knows, however, 
each of them can be
applied  to a specific type of eqautions, to be sure,
but it does not have generality  applicable to various studies.  

                 In this paper we prove a Hopf bifurcation theorem in general Banach spaces,
which improves [\refCRhopf, Theorem 1.11] and can be applicable to semilinear and 
quasi-linear partial differential equations in unbounded domains of $\R^n$. 
Here, we mention some previous results closely reated to our results.
In [\refKiel] Kielh{\"o}fer proved another version of [\refCRhopf, Theorem 1.11] by
using the spaces of 2$\pi$- periodic H{\"o}lder continuous functions which are described
in Section 2 of this paper. This theorem also needs, however, some compactness condition.
So, it can not be applied to partial differential equations in unbounded domains of $\R^n$.
In [\refKhi] the author proved a Hopf bifurcation theorem in Hilbert spaces, which
improves  [\refCRhopf, Theorem 1.11] and can be applicable to semilinear 
partial differential equations in unbounded domains of $\R^n$.
It seems, however, to be difficult to apply the bifurcation theorem to quasi-linear
partial differential equations.

                We consider the next abstract semilinear equation in Banach spaces in this paper:
\[
u_t = A u + h(\lam, u),                \tag{\eqpde}		
\]
where $\lambda$ is a real parameter. The linear operator $A$ and the map  $h$ are described in Section 2 below. 

                     The assumptions of our main theorem ({\thmmain} below) are
weaker than those of   [\refCRhopf, Theorem 1.11].
Actually, our result has the following features:
 \begin{itemize}
  \item 
We do not assume  that  $A$ generates a $C_0$-semigroup.
  \item
We do not assume  that 
 $A$ has compact resolvents.
   \end{itemize}
These features  contribute to  wider applications (see Section 5 below).
In particular, the latter feature makes it possible to apply
our main theorem  to nonlinear partial differential equations on 
unbounded domains of $\R^n$. Actually, we treat the Cauchy problems 
for  semilinear and quasi-lineaer heat systems as 
concrete examples in Section 5 below.

               The idea of the proof of our main theorem in this paper is 
the same as that of the main theorem in [\refKhi].
Actually, the both proofs are based on [\refKsome, Theorem 3].
The technical aspect of our proof in this paper is, however,  
 more complicated.
In [\refKhi] Parseval's identity plays an important role, which does not
hold in general Banach spaces. To overcome the technical difficulty,
we use the H{\"o}lder spaces introduced in [\refKiel] and [\refABB, Theorem 4.2]
which is a result on the well-posedness of linear differential equations in 
H{\"o}lder spaces.

                   The plan of our paper is the following.
In Section 2 we describe our main results and discuss the features of our
results. We describe some preliminary results to prove our main results
in Section 3. We prove our main result in Section 4. In Section 5
we present a concrete example.

\newchapter

\centerline{\bf  2. Hopf bifurcation theorem}

                       Let $V$ be a real Banach space and $V_c= V + i V$ be  its 
complexfication.
Let $A$ be a closed linear operator on $V$ 
with a bounded inverse $A^{-1}$.
We denote its domain by $\domain(A)$, range by $\range(A)$,
null space by $\kernel(A)$ and the extension of $A$ on $V_c$ 
by $\Ac$.  We use the same notation  for the complexfication of
the other linear operaters.
If $W$ is another Banach space,
$\linear(V, W)$ denotes the set of bounded linear operators from
$V$ to $W$. We simply write $\linear(V):= \linear(V, V)$.
We define the real Banach space $U:=\domain(A) \subset V$
with the norm $\| u \|_U := \|  A u \|_{V}$ for $u \in U$.
Let $\beta \in (0,1)$.
We set the real Banach spaces $X$ and $Y$ by
\[
X:= C^{1+\beta}_{2\pi}(\R, V) \cap C^{\beta}_{2\pi}(\R, U)
\moji{and} Y:= C^{\beta}_{2\pi}(\R, V).
                                \tag{\eqholder}
\]
Here, for a Banach space E we denote by $C^{\beta}_{2\pi}(\R, E)$
the space of H{\"o}lder continuous $2\pi$-periodic functions 
$u\semicolon \R \to E$ of 
H{\"o}lder index $\beta$, i.e.
\[
\begin{split}
 C^{\beta}_{2\pi}(\R, E):= & \,                  \biggl\{ u \in C(\R, E) \semicolon 
u(t + 2 \pi) = u(t) \moji{for} t \in \R             \moji{and} \\
  & {\hskip 2.7cm} 
\| u \|_{E, \beta}:= \max_{t \in \R} \| u(t) \|_E + 
                 \sup_{s \not = t} 
                      \dfrac{\| u(t) - u(s) \|_E}{|t - s|^{\beta}} < \infty  \biggr\},   \\
\end{split}
\]
\[
C^{1+\beta}_{2\pi}(\R, E):= \left\{ u \in C^{\beta}_{2\pi}(\R, E) \semicolon 
\dfrac{du}{dt} \in C^{\beta}_{2\pi}(\R, E)  \right\}
\]
with the norm $ \| u \|_{E, 1+\beta}:= \| u \|_{E, \beta} 
                                + \left\|  du/dt \right\|_{E, \beta}
                                \moji{for} u \in C^{1+\beta}_{2\pi}(\R, E)$.
\propspace
                  We assume the following (\HoneFirst) - (\HoneFourth) :
\propspace

\n
{(\HoneFirst)} {\hskip 0.2cm} 
 There exist a real open interval $K$   and $\delta  \in (0,\infty]$
 such that $0 \in K$ and $h$ is a map from $K \times B_U(0 \semicolon \delta)$ to $V$.  
Here, $B_U(0 \semicolon \delta) := \{ u \in U \semicolon \|u\|_U < \delta \}$.
 {\vskip 0.25cm}
\medn
For any $(\lam, u) \in K \times B_X (0 \semicolon \delta)$, we set 
$[h(\lam, u)](t) := h(\lam, u(t)) \in V$ for any $t \in \R$.
\n

 {\vskip 0.25cm}
\medn
(\HoneSecond) {\hskip 0.2cm} 
$h(\lam, u) \in Y$ for any $(\lam, u) \in K \times B_X(0 \semicolon \delta)$.
{\vskip 0.25cm}
\medn
We define the map $\Psi \colon (\lam, u) \in K \times B_X(0 \semicolon \delta)
\mapsto h(\lam, u) \in Y$.
 {\vskip 0.25cm}
\medn
(\HoneThird) {\hskip 0.2cm} 
$\Psi \in C^2(K \times B_X (0 \semicolon \delta), \,Y)$.
 {\vskip 0.25cm}
\propspace
{\bf \remcomm}.
We can regard $U$ (resp. $V$) as the closed subspace of 
$X$ (resp. $Y)$ which consists of constant functions in $X$ (resp.$Y)$.
Then we verify that  (\HoneThird) implies 
$h \in C^2(K \times B_U (0 \semicolon \delta_1), \,V)$ 
for some $\delta_1 > 0$ with
\[
[\Psi_u (\lam, u) v](t)= h_u (\lam, u(t)) v(t),                              
{\hskip 0.4cm}
[\Psi_{uu} (\lam, u) v w](t)= h_{uu} (\lam, u(t)) v(t) w(t)   \moji{in} V                        
\]
and so on for $\lam \in K$, $u, v, w\in X$ and  $t \in \R$.               \End
\propspace
\n
(\HoneFourth) {\hskip 0.2cm} 
$h_u(0, 0)=0$ and $h(\lam, 0) = 0$   if  $\lam \in K$.
\medn
In what follows we simply denote  (\HoneFirst) - (\HoneFourth)  
by  (\Hone). We also assume (\Htwo) - (\Hfive) below.
\medn
(\Htwo)  $\,\,\,\pm\, i$ are the  simple eigenvalues of  $A$, i.e. 
\[
\begin{cases}
& \dim \kernel (i - \Ac) = 1  =\codim \range(i - \Ac), \hfill  \\ 
& \psi \in \kernel(i - \Ac) - \{ 0 \} \,\,
                                \Longrightarrow \,\, \psi \not \in \range(i-\Ac). \hfill \\ 
\end{cases}
\]
So, by the implicit function theorem,
  $\Ac + \{ h_u(\lam, 0) \}_c$ has an eigenvalue $\mu(\lam) \in \C$
and eigenfunction $\psi(\lam) \in \domain(\Ac)$
corresponding to $\mu(\lam)$  
for any $\lam$ in a small  neighborhood of $0$ 
such that $\mu(0)= i$ and that 
$\mu(\lam)$ and $\psi(\lam)$
 are functions of class $C^2$. 
\medn
(\Hthree)\q               (Transversality condition of eigenvalues)\q 
$\text{Re}\,\mu'(0) \not= 0$,
\medn
(\Hfour) \q $i k \in \rho(\Ac)$
                               for  $k\in \Z - \{ -1, 1 \}$.
\medn
 (\Hfive)   There exists $M \in (0, \infty)$ such that
\[
\| (i n - A_c)^{-1} \|_{V_c \to V_c}  \le \dfrac{M}{n} \moji{for} n= 2, 3, 4,\cdots.
\]

                                To begin with, we shortly state our result:
\medn
{\bf  \propmain}.  
{\sl  Let $V$ be a real Banach space and $A$ be a closed linear operator
on $V$.
We assume (\Hone) - (\Hfive). 
Then, $(\lam, u)= (0,0)$ is a Hopf bifurcation point of (\eqpde).
}
\medbreak

             {\propmain} is a short  version of our main result {\thmmain} below, 
which shows that the branch of bifurcating periodic solutions  are
unique in a neighborhood of $(\lam, u)= (0,0)$.

              Next, we make preparation to state our main result.
Let $m \in \Z$, $n \in \N$  and $u \in V_c$.
We write $e_m(t) := e^{imt}$, $c_n(t):= \cos nt$ and
$s_n(t):= \sin nt$  for $t \in \R$.
We denote $(u \otimes e_m)(t) := u e_m(t) = u e^{imt}$ ($t \in \R$).
 Similarly, $(u \otimes c_n)(t) := u \cos nt$
and $(u \otimes s_n)(t) := u \sin nt$ ($t \in \R$).                                                     
We set
$X_1 := \{ u \otimes c_1 + v \otimes s_1 \semicolon u, v \in U \}$
   as a subspace of $X$. 
We define the translation operator $\tau_\theta$ by
$(\tau_\theta u)(t):= u(t - \theta)$ for any $\,\theta \in \R$.
                                                                                                                                                                
               For simplicity, we set $f(\lam, u)= Au + h(\lam, u)$.
If $u(t)$ is a $2\pi$-periodic solution of the next equation (\eqtpe)
then $u(t/(\sigma + 1))$ is a $2\pi (\sigma + 1)$-periodic 
solution of (\eqpde):
\[
u_t = (\sigma + 1) \{ Au + h(\lam, u) \}     .                               \tag{\eqtpe}
\]

                            Our main theorem is the following:
\propspace
{\bf  \thmmain}.  
{\sl We assume (\Hone) - (\Hfive). 
Then, there exist $a, \e > 0$, $u_{\star} \in X_1 - \{ 0 \}$ 
and functions
$\zeta = (\lam, \sigma) \in C^1([0, a), \R^2)$, 
$\eta \in C^1([0, a), X)$ with the following properties:
\medn
(a)  \,\, $(\lam, \sigma, u)= (\zeta(\alpha), \alpha u_{\star} + \alpha \eta(\alpha) )$
is a $2\pi$-periodic solution of (\eqtpe),
\medn
(b) \,\, $\zeta(0)=\zeta'(0)=(0, 0)$ and $\eta(0)=0$,
\medn
(c) \,\, If $(\lam, v)$ is a solution of (\eqpde) of period 
$2\pi(\sigma + 1)$, $|\lam| < \e$, $|\sigma| < \e$,
${\tilde v} \in X$ and
$\|  {\tilde v} \|_X < \e$, where
${\tilde v}(t):= v((\sigma + 1)t )$ for $t \in \R$, then there exist
$\alpha \in (0, a)$ and $\theta \in [0, 2\pi)$ such that 
$(\lam, \sigma)= \zeta(\alpha)$ and $\tau_\theta {\tilde v}= \alpha u_{\star}
+ \alpha \eta(\alpha)$.
}

\npage

\centerline{\bf 3. Preliminary results}

\medbreak
First, we describe a basic bifurcation theorem ({\thmbasic} below), which is a slightly refined 
version of [\refKsome, Theorem 3]
for the case $m=2$ . The proof of our main result (\thmmain) is based on 
this result.

Let $\cX$ and  $\cY$ be real Banach spaces and 
$O$ be an open neighborhood of $0$ in $\cX$.
Let $J$ be an open neighborhood of $0$ in $\R^2$.
Let  $g \in C^2( J \times O, \,\cY)$ be a 
map such that
\[
g(\Lam, 0)= 0 \moji{for \,\, any}               \Lam= (\Lam_1, \Lam_2) \in J.
\]
We define $H \colon J \times \cX \to \R^2 \times \cY$ by
\[
H
\begin{pmatrix}
\Lam  \\
u  \\
\end{pmatrix}
:=
\begin{pmatrix}
l u - \eone  \\
g_u(\Lam, 0)u  \\
\end{pmatrix}.                   \tag{\eqH}
\]
Here, $l :=(l^1, l^2) \in \linear(\cX, \R^2)$
and $\eone :=  (1, 0)$.
We define 
$G \colon J \times O\, \to \,     \R\times \cY$ by
\[
G
\begin{pmatrix}
\Lam  \\
u  \\
\end{pmatrix}
:=
\begin{pmatrix}
l ^2 u  \\
g(\Lam, u)  \\
\end{pmatrix} .
\] 
We set  $Z:= \kernel(l)=  \{ u \in \cX \semicolon l u= (0, 0) \}$.

\propspace
{\bf \thmbasic}. 
{\sl In addition to the 
assumptions above we assume that $\us \in O$ satisfies
\medn
{\rm (\eqThe)}\q  $(\Lam, u)= (0,\us)$ is an isolated solution of the
extended system $H(\Lam, u) = 0$. 
\medn
 Then there exist an open neighborhood $W$ of (0,0)  
 in $\R^2 \times \cX$,
$a  \in (0, \infty)$ and  functions 
$\zeta \in C^1((-a, a), \R^2)$, 
$\eta \in C^1((-a, a), Z)$ such that
 $\zeta(0)=0$, $\eta(0)=0$ and
\[
  G^{-1} (0) \cap W   = \{ (\Lam, 0) \semicolon  (\Lam, 0) \in W \}                           
\cup \{  (\zeta(\alpha), \, \alpha \us + \alpha\eta(\alpha)) \semicolon |\alpha| < a \} . 
                                                                                                                                                                                                                                                                                                 \tag{\eqgi}
\]
}
  \medbreak
          {\it Proof}. We set ${\tilde Z} := \kernel(l^1)$.
By [\refKsome, Theorem 3] the statement of {\thmbasic}
with $\eta \in C^1((-a, a), Z)$ replaced by
 $\eta \in C^1((-a, a), {\tilde Z})$ holds. It follows that 
$G(\zeta(\alpha), \alpha \us + \alpha \eta(\alpha)) = 0\,$
 for any $\alpha \in (-a, a)$, which implies
$\,0= l^2\{ \alpha \us + \alpha \eta(\alpha) \} = 
\alpha l^2 \eta(\alpha)$.                          
Therefore, $\eta(\alpha) \in Z$ for any $\alpha \in (-a, a)$.
                                \End
\eqspace

    Next, we use the same notation in Section 2. We set
$Y_1 := \{ u \otimes c_1 + v \otimes s_1 \semicolon u, v \in V \}$
   as a subspace of $Y$. 
   We define 
$L_1 : \Vc \to Y_1$ by
$L_1\psi := \re (\psi \otimes e_1)$ 
for any $\psi \in \Vc$
 and 
$T_1 : X_1 \to Y_1$  by
$T_1 w := \dot{w} - Aw$ for any 
$w \in X_1$. 
    Then, it follows that 
\begin{gather*}
L_1(a + i b) = a \otimes c_1 - b \otimes s_1
\moji{for any} a, b \in V,            \tag{\eqlon}        \\ 
T_1 (a \otimes c_1 + b \otimes s_1)
= (b - A a)\otimes c_1 - (a + A b)\otimes s_1
                                                                \moji{for any} a, b \in U.            \tag{\eqtoa}
\end{gather*}
 In view of (\eqlon) the following result clearly holds:
\medn
{\bf \propisom}. {\sl  If we regard   $\Vc$ and $\Uc$ as  real linear spaces then
we have the following results.
{\vskip 0.2cm}
\n
(i) The operator $L_1$ is  isomorphic as a real linear operator
from  $V_c$ to  $Y_1$. 

{\vskip 0.2cm}
\n
(ii) The operator
 $L_1 |_{\Uc} $ is  isomorphic as a real linear operator
 from  $U_c$ to  $X_1$. 
 }

\propspace

\n
{\bf \propkernel}. {\sl  
(i) \, $L_1\, \kernel (i - \Ac) =  \kernel(T_1)$,
\medn
(ii) \, $L_1 \range (i - \Ac) =  \range(T_1)$.
}
\medbreak
                {\it Proof.} 
If $w \in X_1$, by  {\propisom} (ii) there exists a unique $\psi\in U_c$
such that $w = L_1 \psi$. Then, we  verify that $T_1 w = L_1(i - \Ac)\psi$,
which clearly leads to (i) and (ii).        \End
\vspace{0.2cm}                                                                                                                     
\eqspace
{\bf \proppsd}.  {\sl  Let $\psi \in \Uc$ and $w= L_1\psi$.
 \medn
 (i) {\hskip 0.2cm} $L_1 (i \psi)= \dot{w}$,
 \medn
 (ii) If $\psi \in \kernel(i - \Ac)$,  then 
 $\,L_1 (i \psi)= A w$.
 }
  \medbreak
                                                                                                {\it Proof}. (i)  
 $L_1(i \psi)= \text{Re\,}\biggl[  \dfrac{d}{dt} 
 (\psi  \otimes e_1 )\biggr]
 = \dfrac{d}{dt} \,\text{Re\,} (\psi  \otimes e_1) = \dot{w}$.
\par
            (ii) We immediately obtain the desired conclusion from (i)
and {\propkernel} (i).                               \End

\propspace\par

         Finally, we assume (\Hone) - (\Hfive).
let $\psis \in \kernel(i - A_c) - \{ 0 \}$ (see (\Htwo)).
We set $\Vst:= \spnc \{ \psis, \overline{\psis} \}$
and  $\Vsh:= \range(i - \Ac) \cap \range(- i - \Ac)$,
which are closed subspaces of $\Vc$
 (see e.g. [\refSpectral, Theorem 3.2]
 for the closedness of $\Vsh$). 
 Let $\Ast := \Ac|_{\Vst}: \Vst \to \Vst$ and
$\Ash := \Ac|_{\Vsh}: \Vsh \to \Vsh$.
Here, $\Ash$ is well-defined by the following lemma:
\propspace
{\bf  \lemwel}.  $\Ac \, \{ \Vsh \cap\Uc \} \subset  \Vsh$.
\propspace

            {\it Proof}.  Let $\varphi \in \Vsh \cap \Uc$. Then,  we have
$\Ac \varphi = i \varphi - (i - \Ac)\varphi \in \range(i - \Ac)$
and
$\Ac \varphi = - i \varphi - ( - i - \Ac)\varphi \in \range(- i - \Ac)$.
So, $\Ac \varphi \in \Vsh$.         \End
                
\propspace
We verify that $A_\star$ and $A_\sharp$ are closed operators.  
We have the following results:  

\propspace
{\bf  \propvzc}.  
{\sl 
(i) $\Vc = \Vst \directsum \Vsh$.
\medn
(ii) $\Ac = \Ast \directsum \Ash$.
\medn
(iii) $i \Z \subset \rho(\Ash)$.
}

\propspace

          {\it Proof}. (i)    In view of (\Htwo),
\[
\range(i - \Ac)          \directsum \spnc\{ \psi_{\star} \} 
                                = \Vc.                                                       \tag{\eqria}
\]
 Taking the complex conjugate of (\eqria), we have
\[
\range    (- i - \Ac)      \directsum \spnc\{ \overline{\psis} \}        
                                                                  = \Vc.                        \tag{\eqrsc}
\]
Let $\varphi \in \Vc$. In view of (\eqria) and (\eqrsc), 
there exist $c_1, c_2 \in \C$,  $\varphi_1 \in \range(i - \Ac)$
and  $\varphi_2 \in \range( - i - \Ac)$  such that
\[
\varphi = c_1 \psis + \varphi_1 
                   \moji{and}                \varphi_1 = c_2 \overline{\psis} + \varphi_2.                           
                                                                                                  \tag{\eqpco}
\]
 It follows that
\[
\varphi = c_1 \psis + c_2 \overline{\psis} + \varphi_2.             \tag{\eqvarp}
\]
By the second equality of (\eqpco) and 
$\overline{\psis} = (1/ 2 i)(i - \Ac)\overline{\psis} \in \range(i - \Ac)$,
we have
$\varphi_2 \in\Vsh$.
So, $\varphi = c_1 \psis + c_2  \overline{\psis} + \varphi_2 \in \Vst + \Vsh$.
It follows that $\Vc = \Vst + \Vsh$. Next, let $\gamma, \delta \in \C$
and $\psish \in \Vsh$ satisfy 
$\gamma \psis + \delta \overline{\psis} + \psish = 0$.
It follows from (\eqria) that
 $\gamma = 0$ and $\delta \overline{\psis} + \psish = 0$.
 By (\eqrsc),  we have $\delta = 0$ and $\psish=0$.
 Thus, $\Vc = \Vst \directsum \Vsh$.
\par
                    (ii)   In view of (i), we define the projection
$P \in \linear(\Vc)$ onto $\Vst$.
Then,  we verify that $P \Ac \subset \Ac P$. 
So, we have the desired conclusion. 
\par
                   (iii)   In view of (ii) and (\Hfour),  
$i k \in \rho(\Ash)$ for $k \in \Z - \{ -1,1 \}$.
So, it suffices to show $i \in \rho(\Ash)$.
By (i) and (\Htwo), $i - \Ash$ is one to one.
Next, let $v \in \Vsh$.  By (i) and $v \in \range(i - \Ac)$
there exist $\alpha, \beta \in \C$ and 
$u \in \Vsh \cap \Uc$ such that
$(i - \Ac)(\alpha \psis + \beta \overline{\psis} + u) = v$.
It follows from (ii) that
$2i \beta \overline{\psis} + (i - \Ash) u = v$ and $(i - \Ash)u \in \Vsh$.
Again by (i), we have $\beta=0$ and $(i - \Ash)u = v$.
So, $i - \Ash$ is  onto and  $i \in \rho(\Ash)$.
                                \End

\propspace
                                                                                                
                              Finally,  let $\cV$ be a complex Banach space and
$\cA$ be a closed operator on $\cV$. 
Let $f \in C^{\beta}_{2\pi}(\R, \cV)$.
 We consider the next problem:
\[
u_{t}= \cA u + f(t) \moji{for} t \in \R.        \tag{\eqPper}
\]
Let $\beta\in (0, 1)$.                                
We have the following result:
\propspace\n
{\bf \propab}.           
{\sl
We assume 
$i \Z \subset \rho(\cA)$ and
$\sup\,\{ \| n (i n - \cA)^{-1} \| \semicolon n \in \Z \} < \infty$.
Then, for each $2\pi$-periodic function $f \in C^{\beta}_{2\pi}(\R, \cV)$
there exists a unique  periodic solution \linebreak
$u \in C^{1+ \beta}_{2\pi}(\R, \cV) \cap C^{\beta}_{2\pi}(\R, \domain(\cA))$ of (\eqPper). 
}
\propspace

                 {\it Proof}. We denote by $C^{\beta}_\text{per}([0, 2\pi], \cV)$
the space of H{\"o}lder continuous  functions $u\colon [0, 2\pi] \to \cV$ of 
H{\"o}lder index $\beta$ such that $u(0) = u(2\pi)$.
We consider the next problem:
 \[
\begin{cases}        
& {\hskip - 0.25cm}  u_{t}= \cA u  + f(t) \moji{for} t \in [0, 2 \pi],   \\ 
& {\hskip - 0.25cm}       u(0) = u(2\pi). \\ 
\end{cases}
             \tag{\eqEq}
\]
By  [\refABB, Theorem 4.2] there exists a unique periodic function
$u= {\tilde u} \in C^{1+\beta}_\text{per}([0, 2\pi], \cV)$
        $\cap\,\, C^{\beta}_\text{per}([0, 2\pi],  \domain(\cA))$.
We denote by $\hat{u}(t)$ 
the $2\pi$-periodic extension of
$\tilde{u}$ to $\R$. Then, $u(0)=u(2\pi)$ implies that $u_t(0)=u_t(2\pi)$. 
  So, $\hat{u} \in C^{1+\beta}_{2\pi}(\R, \cV) 
                \,\cap \, C^{\beta}_{2\pi}(\R, \domain(\cA))$, which is a unique 
solution of  (\eqPper).      \End

\newchapter
 
 \centerline{\bf 4. Proof of {\thmmain}}
\medbreak 
                            Let $X$ and $Y$ be  real Banach spaces defined by (\eqholder).
We denote the $n$-th Fourier coefficient  of $\varphi \in Y_c$ by
\[
{\hat \varphi}(n):= \dfrac{1}{2\pi}\int_{0}^{2\pi}  \varphi(t) e^{-int} dt. \q
        (n \in \Z)           \tag{\eqfrac}
\]
We set 
\[
X_0 :=U 
\moji{and} X_{\infty} := 
\{ \varphi \in X \semicolon {\hat \varphi}(n)=0 
\moji{for} n=-1, 0, 1 \}                                                              \tag{\eqxn}
\]
  as closed  subspaces  of $X$,
\begin{align*}
Y_0 :=V, \,\,\,\, &
Y_1 := \{ u \otimes c_1 + v \otimes s_1 \semicolon u, v \in V \}         \tag{\eqyn}  \\ 
&  \moji{and} Y_{\infty} := 
\{ \varphi \in Y \semicolon {\hat \varphi}(n)=0 
\moji{for} n=-1, 0, 1 \}  
\end{align*}
as  closed subspaces of $Y$.
Let $X_1$  be a closed subspace of $X$ 
defined in Section 2.

\medbreak
     {\it Proof of {\thmmain}}.               We apply {\thmbasic}.
We use the notation in Section 2 and 3.
We denote $\Lam= (\lam, \sigma) \in \R^2$.
We define $g \in C^2(K\times\R \times B_X(0, \delta), \,\,Y)$ by
$g(\Lam, u)= u_t - (\sigma +1) f(\lam, u)$,
where $f(\lam, u):= A u + h(\lam, u)$.
By the assumption (\Htwo) in Section 2 there exists  
  $\,\psi_{\star} \in \kernel(i - \Ac) - \{ 0 \}$.
Then, $\re \psis$ and $\im \psis$ are linearly
independent in $V$.
 So, by the Hahn-Banach theorem
there exists $\,m \in V^*\,$ such that
 $m_c \psi_{\star} = 1$. 
We define  $l = (l^1, l^2) \in \linear(X, \R^2)$ by
\[
l^1 u  := \frac{1}{\pi}\int_0^{2\pi} m u(t) \cos t \,dt
\,\, \moji{and}           \,\,
        l^2 u  := \frac{1}{\pi}\int_0^{2\pi} m u(t)  \sin t \,dt
\]
for $u \in X$.
We set
$ u_\star := L_1 \psi_{\star} =   \re(\psis \otimes e_1) \in X_1$.               
Then, 
 $l u_{\star} = (1,0) = {\bm e}_1$.                                                       
Let $H \colon K \times \R \times X \to \R^2 \times Y$ be the  operator defined by (\eqH).
Then, by  (\Hone-4) and {\propkernel} (i),
$H(0, u_{\star}) = 
(l u_{\star} - {\bm e}_1, \, (\us)_t - A\us) = (0,0)$.
We set $DH^{\star}:= DH(0, \us)$.  
Then, we have
\[
DH^\star
\begin{pmatrix}
\lam  \\ 
\sigma  \\ 
u  \\ 
\end{pmatrix}
=
\begin{pmatrix}
l^1 u \\ 
l^2 u \\ 
u_t - Au - \sigma A \us - \lam h_{\lam u}^0 \us  \\ 
\end{pmatrix} ,         \tag{\eqmapdh}
\]
where $h_{\lam u}^0:= h_{\lam u}(0,0)$.
We verify that 
$S:= DH^\star|_{\R^2 \directsum X_0 \directsum X_1} :
\R^2 \directsum X_0 \directsum X_1 \to 
\R^2 \directsum Y_0 \directsum Y_1$ 
and 
$T:= DH^\star |_{X_{\infty} } :
X_{\infty} \to Y_{\infty}$ 
are well-defined by (\eqmapdh)
and that $DH^\star = S \directsum T$.
We note that
$T u = u_t - A u$ for any $u \in X_{\infty}$.
In view of  the below {\lemsi}
and {\lemti}, $DH^{\star}$ is bijective.
So,  by {\thmbasic} 
 $(\lam, u)= (0, 0)$
is a Hopf bifurcation point of (\eqpde) and 
 there exist an open neighborhood $W$ of (0,0)  
 in $\R^2 \times X$,
$a  \in (0, \infty)$ and  functions 
$\zeta \in C^1((-a, a), \R^2)$, 
$\eta \in C^1((-a, a), Z)$ such that $\zeta(0)=0$, $\eta(0)=0$ and
(\eqgi) holds. 
Here, $Z:= \{ u \in X \semicolon l u=(0,0)  \}$. 
So, (a) holds. 
Next, we show the following (4.5) in preparation
to prove (b) and (c).
\[
\zeta(- \alpha) = \zeta(\alpha) 
\moji{and}
 \eta(- \alpha) = - \tau_\pi (\eta(\alpha)) 
                  \moji{for any} \alpha \in [0, a).           \tag{\eqgam}
\]
We set $U(\alpha):= \alpha \us + \alpha \eta(\alpha) \in X$ for any 
$\alpha \in (-a, a)$.
We define $V(\alpha) \in X$ by $V(\alpha):= \tau_\pi (U(\alpha))$.
Let  $\gamma \in (0, a)$ be a constant such that 
$\{ (\zeta(\alpha), V(\alpha)) \semicolon \alpha \in [0, \gamma) \} \subset W$.
Then, $(\zeta(\alpha), V(\alpha)) \in G^{-1}(0) \cap W$
for any $\alpha \in [0, \gamma)$. 
So,  by {\thmbasic} for any $\alpha \in [0, \gamma)$
 there exists $\beta \in (-a, a)$ such that
  $(\zeta(\alpha), V(\alpha)) = (\zeta(\beta), U(\beta))$.
  On the other hand, $l^1 V(\alpha) = - \alpha$ and $l^1 U(\beta) = \beta$.
  Therefore, $\beta = - \alpha$ and $(\zeta(- \alpha), U(- \alpha)) = (\zeta(\alpha), V(\alpha))$
  for any $\alpha \in [0, \gamma)$. 
  Actually, we easily verify from the frequently used argument by 
  contradiction that 
\[
a = \sup\,\{ q \in (0, a) \semicolon (\zeta(- \alpha), U(- \alpha)) 
    = (\zeta(\alpha), V(\alpha)) 
                                \moji{for any} \alpha \in [0,q) \}.    \tag{\eqdefa}
\]
We obtain (\eqgam) from (\eqdefa) and $\tau_\pi \us = - \us$. 
 
            By  (\eqgam), $\zeta'(0)=(0, 0)$. So, (b) holds.
                              Finally, we show (c).
Let $\e$ be a positive constant
such that if $(\lam, \sigma, w) \in \R^2 \times X$ 
satisfies $|\lam| < \e, |\sigma| < \e$
and $\| w \|_{X} < \e$ then $(\lam, \sigma, w) \in W$.
Now, let $(\lam, v)$ be a solution of (\eqpde) of period 
$2\pi(\sigma + 1)$, $|\lam| < \e$, $|\sigma| < \e$,
${\tilde v} \in X$ and
$\|  {\tilde v} \|_X < \e$,  where
${\tilde v}(t):= v((\sigma + 1)t )$ for $t \in \R$.
For simplicity, we set $(p, q):= l {\tilde v}=(l^1 {\tilde v}, l^2 {\tilde v})$.
First we consider the case: $q = 0$.
Then $(\lam, \sigma, {\tilde v}) \in W$ is a  solution of 
$G(\Lam, u):=(l^2 u, g(\Lam, u)) = (0, 0)$.
By {\thmbasic} there exists $\alpha \in (-a, a)$ such that 
$(\lam, \sigma) = \zeta(\alpha)$ 
and ${\tilde v} = \alpha \us + \alpha \eta(\alpha)$.
If $\alpha < 0$ then 
$(\lam, \sigma) = \zeta(-\alpha)$ and 
$\tau_\pi {\tilde v} =   (- \alpha) \us + (- \alpha) \eta(- \alpha)$
in view of (\eqgam) and $\tau_\pi \us = - \us$.
           Next, we consider the case: $q \not = 0$.
There exists   $\theta \in (0, 2\pi)$ such that 
$e^{i \theta}=(p - i q)/\sqrt{p^2+ q^2}$.
Then, $l^2 \tau_\theta {\tilde v} = 0$ and 
$(\lam, \sigma, \tau_\theta {\tilde v}) \in W$ 
is a solution of $G(\Lam, u)=0$. 
 So, the present case is reduced to the case: $q = 0$.
Therefore, (c) holds.                                \End           
 \propspace
                    In the above proof, we use the following two lemmas:
 \eqspace
{\bf \lemsi}. {\sl  The operator $S$ is bijective.}
\eqspace
{\bf \lemti}. {\sl  The operator $T$ is bijective.}
\eqspace 
             {\it Proof of {\lemsi}}.  
The idea of proof is essentialy the same as that of [\refKhi, Lemma 4.1].
By  (\Htwo), {\remcomm} and the implicit function theorem
(see e.g. [{\refCRbif}, Theorem A])
$\{ f_u(\lam, 0) \}_c$ has an eigenvalue $\mu(\lam) \in \C$
and an eigenfunction $\psi(\lam) \in U_c$ 
corresponding to $\mu(\lam)$ 
for any $\lam$ in a small open interval $K_1$ 
such that  
$0 \in K_1 \subset K$,  $\mu(0)= i$, $\psi(0)= \psi_{\star}$,
$\mu(\cdot) \in C^2(K_1, \C)$ and $\psi(\cdot) \in C^2(K_1, U_c)$.
Differentiating 
$\{ f_u(\lam, 0) \}_c \, \psi(\lam)=\mu(\lam)\psi(\lam)$
with respect to $\lam$, we have
\[
\mu'(0)\psi_{\star} + (i  - \Ac) \psi'(0)                   
             = (f^0_{\lam u})_c  \, \psi_{\star}.     \tag{\eqmuz}
\]
We set $p:= \re \mu'(0)$ ($\not = 0\,$  by (\Hthree)), 
$q= \im\mu'(0)$ and
$u_{\sharp}:= L_1 \psi'(0) \in X_1$.
 It follows from (\eqmuz) and 
{\proppsd} that
{\hskip 0.2cm}
\[
 f_{\lam u}^0 \us = p \us + q A \us + T_1 u_{\sharp}.                \tag{\eqfly}
\]
Let $u_0 \in X_0$,  $u_1 \in X_1$ and $u= u_0 + u_1$. 
In view of  (\eqmapdh) and (\eqfly), we have
\[
S                              
\begin{pmatrix}
 (\lam, \sigma) \\ 
u_0    \\
u_1      \\
\end{pmatrix}                                           
=
\begin{pmatrix}
 l u_1 \\ 
 - A u_0   \\
T_1(u_1 - \lam u_{\sharp}) - \lam p \us 
  - ( \sigma + \lam q) A \us      \\
\end{pmatrix}
.                                                               \tag{\eqsvec}
\]
By (\Htwo),  we have
$\range(i - \Ac) \directsum \spn \{ \psis \} =  \Vc$.
It follows from  {\propisom} (i), {\propkernel} (ii) and {\proppsd} 
that
\[
\range(T_1)  \directsum \spn\{ \us, A\us \}  = Y_1.       \tag{\eqrto}
\]
We note that $\us$ and  $A\us$ are linearly independent in $Y$.

                 First, we show that $S$ is one to one.
Let $S(\lam, \sigma, u)= 0$.
It follows from (\Hthree),  (\eqsvec), (\eqrto) and $0 \in \rho(A)$
that $u_0=0$, $\lam= \sigma=0$,
\[
l u_1=(0,0) \moji{and} T_1 u_1=0.                \tag{\eqluo}
\]
Let $\psi_1 := L_1^{-1} u_1 \in U_c$. 
Then by (\eqluo)  and {\propkernel} (i),
\[
\psi_1 \in \kernel(i - \Ac) \moji{and} m_c \psi_1 = 0.                                                
                                                                                                                                                                                                                                                                                                                                                                                                 \tag{\eqpin}
\]
It follows from (\eqpin), (\Htwo) and $ m_c \psis  = 1$  that 
$\psi_1 = 0$, which implies $u_1=0$. 
So, $S$ is one to one. 
\par
                            Next, we show that $S$ is onto.
Let $(a, b, y_0, y_1) \in \R^2 \directsum Y_0 \directsum Y_1$.
In view of  $\,0 \in \rho(A)$,
there exists $x_0 \in X_0$ such that $- A x_0 = y_0$.
By (\eqrto) there exist $w \in \range(T_1)$ and 
$(\gamma, \delta) \in \R^2$ such that 
\[
w + \gamma u_{\star} + \delta A u_{\star} = y_1.                                                     \tag{\eqwg}
\]
We set $\lam_0:= - \gamma/p$ and 
$\sigma_0 := -\delta + \gamma q /p$.
There exists $v_1 \in X_1$ such that 
$T_1 (v_1 - \lam_0 u_{\sharp}) = w$.
Let $(\alpha, \beta) := l v_1 \in \R^2$ and
$x_1:= v_1 + (a-\alpha)u_{\star}+(\beta - b)Au_{\star} \in X_1$.
By {\propkernel} (i) and {\proppsd} (ii),
we have 
$A u_{\star} = L_1 (i \psi_{\star}) \in \kernel(T_1)$.
So, $lAu_{\star}=(0, -1)$.
It follows from $l \us = \eone$, {\propkernel} (i),
(\eqsvec) and (\eqwg) that
$S(\lam_0, \sigma_0, x_0, x_1)=(a,b,y_0,y_1)$.
Therefore, $S$ is onto.                            \End
\propspace

                  {\it Proof of {\lemti}}.  \q
 It suffices to show that $T_c \colon X_{\infty c} \to Y_{\infty c}$ is bijective.
                                                                Let $v \in Y_{\infty c}$.
 We will show that the following equation (\eqtuv) has a unique solution 
 $u \in X_{\infty c}$.
\[
T_c u = v \q (\Longleftrightarrow u_t - A_c u = v)            \tag{\eqtuv}
\]

            To begin with, we consider the uniqueness of solutions for (\eqtuv).
let $v=0$.
 By the Fourier transform we have
$(i n - A_c) {\hat u}(n) = 0 \moji{for any} n \in \Z$.
In view of (\Hfour) and $u \in X_{\infty c}$, ${\hat u}(n) = 0$ 
for $n \in \Z$. 
So, $u=0$ by Fejer's theorem, which implies
the uniqueness of solutions for (\eqtuv). 

                  Next, we consider the existence of solutions for (\eqtuv).
   Let $v \in Y_{\infty c}$.
In view of  {\propvzc} (i) we can define the projection 
$P\in \linear(V_c)$ onto $\Vst$.
We decompose (\eqtuv) into the following two equations:
\[
u_t - \Ast u = P v \moji{on} \Vst,                                 \tag{\eqwta}
\]
\[
u_t   - \Ash u = (I - P)v  \moji{on} \Vsh.                  \tag{\eqzta}
\]

              First, we consider (\eqwta). 
There exist $g, h \in C^{\beta}_{2\pi}(\R,\, \C)$ such that
 $P v(t) = g(t) \psis +
h(t) \,\overline{\psis}$ for any $t \in \R$.                                  
We set $u_1(t)  := c(t) \psis + d(t) \, \overline{\psis} \in X_{\infty c}$.
Then, by substituting  $u=u_1(t)$ for (\eqwta) we have
\[
c\,'(t) - i c(t) = g(t) \moji{and}  d\,'(t) + i d(t) = h(t).       \tag{\eqcd}
\]
Considering the condition $\hat{c}(1)= \hat{d}(-1)=0$ we solve 
(\eqcd) to obtain 
\[
c(t):= e^{i t} \{ \varphi_1(t) - \hat {\varphi}_1(0)   \} \moji{and} 
d(t):= e^{- i t} \{  \varphi_2(t) - \hat{\varphi}_2(0)  \}.                         \tag{\eqcdt}
\]
Here, we set $\varphi_1(t) := \int_0^t g(s) e^{-is} ds$, 
$\varphi_2(t) := \int_0^t h(s) e^{is} ds$, which are $2\pi$-periodic fucntions.
We verify that $u_1 \in X_{\infty c}$ and that $u = u_1(t)$ is actually 
a solution of (\eqwta).                                                            

               Next, we consider (\eqzta). 
By {\propvzc} (iii) and (\Hfive), $i \Z \subset \rho(\Ash)$ and \linebreak
   $\sup_{n \in \Z} |n| \cdot 
 \| (i n - \Ash)^{-1} \|_{\Vsh \to \Vsh} < \infty$. 
So, by  {\propab}, the equation
(\eqzta) has a  solution 
$u= u_2 \in C^{1 + \beta}_{2\pi}(\R,V_\sharp)  \cap C^{\beta}_{2\pi}(\R, \domain(\Ash))$.
 In view of  {\propvzc}  (iii) we verify that $\hat{u}_2(n)=0$ for $n=0, \pm 1$. 

               Therefore, $u= u_1 + u_2 \in X_{\infty c}$ is a solution of (\eqtuv).                        
\End

\newchapter

\centerline{\bf 5. A concrete example}

\propspace

                       In this section we freely use the notation used in Section 4.	
\propspace

                     We consider the Cauchy problem of  the following  quasi-linear system:
\[
\begin{cases} 
& {\hskip - 0.3cm} u_{t}=   \{\kappa_1(u) u_x \} _x  - v
	- p u + u (\lam\,q^2 - u^2 - v^2) 
	         \moji{for} (x, t) \in \R \times [0, \infty),   \hfill \\ 
&  {\hskip - 0.3cm} v_{t}= \{\kappa_2(v) v_x \} _x + u
	- p v + v  (\lam\,q^2 - u^2 - v^2) 
	       \moji{for} (x, t) \in \R \times [0, \infty).   \hfill \\ 
\end{cases}
			\tag{\eqhc}
\]
Here, $p$ and $q$ are functions on $\R$ defined by 
$p(x) := \left\{2 \tanh^2(x/2) - 1  \right\} / 4$
and $q (x):= \sech(x/2)$.
Let $I$ be a real open interval such that $0 \in I$.
We assume that the functions $\kappa_1$ and $\kappa_2$ satisfy
the following conditons:
\medn
(A-1) $\kappa_1$, $\kappa_2 \in C^5(I,\R)$,
\medn
(A-2) $\kappa_1(0) =  \kappa_2(0) = 1$,
\medn
(A-3) $\kappa_1(r) > 0$ and  $\kappa_2(r) > 0$ for any $r \in I$.
\medn
We simply denote (A-1) - (A3) by (A). 
In this section we prove the next result by formulating (5.1) in the form of (\eqpde).
\propspace
{\bf \propqua}.
 {\sl
We assume (A). 
Then  $(\lam, u)= (0,0)$ is a Hopf bifurcation point of (\eqhc).
}
\propspace
{\bf \remqua}. As  preliminary study we consider the case 
where (\eqhc) is semilinear, i.e. the case $\kappa_j(r) \equiv 1$ ($j = 1, 2$). 
In this case, as discussed in [\refKhi, Section 5],
 the branch of periodic solutions of  (\eqhc)
 $\,(u, v)={  (u_{\lam}, v_{\lam})}$ $(\lam > 0)$
bifurcates  at $\lam=0$ from the branch of trivial solutions.
Here, ${ u_{\lam}(x, t) :=  \sqrt{\lam} \,q(x) \cos t }$ and ${  v_{\lam}(x, t) :=}$ 
 ${\sqrt{\lam} \,q(x) \sin t \,}$. Interestingly, in  both of quasi-linear equation (\eqhc)
and the semilinear equation  (\eqhc) with $\kappa_j(r) \equiv 1$ $(j= 1,2)$
 the Hopf bifurcation occurs at the same value $\lambda = 0$.   \End
 \propspace
          We make preparations to prove {\propqua}.

\propspace

            We set $V := L^2(\R) \times L^2(\R)$ and $U :=  H^2(\R) \times H^2(\R)$.
Let $\ubold=(u, v)$.  We define $A \colon V \to V$ by 
\[
A \ubold :=
\begin{pmatrix}
u_{xx}  - v - p u  \\ 
v_{xx}  + u - p v  \\ 
\end{pmatrix}
\moji{for} \ubold \in \domain(A):= U.
\]
Let $d > 0$ is a constant satisfying $(- d, d) \in I$. Here, $I$ is the interval described
in the above condition (A).
We  define  $H \colon B_U(0 \semicolon d) \to V$ and 
$h_0, h$ : $\R \times B_U(0 \semicolon d) \to V$ by
\begin{align*}
&H(\ubold):= 
\begin{pmatrix}
\{ (\kappa_1(u) - 1) u_x \} _x \\ 
\{ (\kappa_2(v) - 1) v_x \} _x
\end{pmatrix}
=
\begin{pmatrix}
\kappa_1' (u) u_x^2 +  \{\kappa_1(u) - 1\} u_{xx}   \\ 
\kappa_2' (v) v_x^2 +  (\kappa_2(v) - 1) v_{xx} 
\end{pmatrix}
,      \tag{\eqhu}        \\
&h_0(\lambda, \ubold) := 
\begin{pmatrix}
u (\lam\,q ^2 - u^2 - v^2) \hfill \\
v (\lam\,q ^2  - u^2 - v^2) \hfill 
\end{pmatrix},
\q	
h(\lambda, \ubold) := H(\ubold) \,\, +  \,\, h_0(\lambda, \ubold)  
\end{align*}
for  $\lam\in \R, \,\, \ubold \in B_U(0 \semicolon d)$.
The above maps are well-defined in view of {\lemhone} below.
So, (\Hone-1) in Section 2 holds. 
In view of (\eqhu) we can formulate (\eqhc) in the form of (\eqpde).
           We can  not  apply [\refCRhopf, Theorem 1.11]   to (\eqhc) since 
the linear operator $A$ does not have compact resolvents.
On the other hand, we will apply our {\propmain} to (\eqhc) to prove {\propqua}.
\propspace

        {\it Proof of \propqua}. 
In view of  {\propmain}  it suffices to  show that  (\Hone) - (\Hfive) hold.
 To begin with, we note that we verified (\Hone-1) by the above discussion and 
 that  (\Htwo) -  (\Hfive) have been verified
in [\refKhi]. Actually, (\Htwo) is the same as [\refKhi, (\Bone)] which was 
verified in [\refKhi, Section 5]. In the same way, 
 (\Hthree)  is the same as [\refKhi, (\Btwo)], 
  (\Hfour) as [\refKhi, (\Bthree)] and
 (\Hfive) as [\refKhi, (\Kone)]. 
We verify (\Hone-2), (\Hone-3) and (\Hone-4) by Lemmas 5.4, 5.5 and 5.6.      \End

\propspace

	We need the below Lemmas 5.1 - 5.3 to prove Lemmas 5.4 - 5.6. 

\propspace
{\bf \lemhone}. {\sl  Let $u \in H^1(\R)$. Then $u \in C(\R) \cap L^\infty(\R)$
with the estimate
\[
\| u \|_{L^\infty(\R)} \le \frac{1}{\sqrt{2}} \, \| u \|_{H^1(\R)}.         \tag{\eqho}
\]
}
 \medbreak
          {\it Proof}. 
Though the inequality  (\eqho) seems to be known,  we sketch its proof
for the sake of completeness. 
              Let $u \in H^1(\R)$. Then   $u \in C(\R) \cap L^\infty(\R)$ by
Sobolev embedding theorem. Since $C_0^\infty(\R)$ is dense in $H^1(\R)$, 
it is sufficent to show (\eqho) under the assumption: $u \in C_0^\infty(\R)$.
Let $a$, $b \in \R$ satisfy that $a < b$ and 
support($u$) is included in the  interval $(a, b)$.
Let $\xi \in (a, b)$. Then,
\[
\{ u(\xi) \}^2 = 2 \int_a^\xi u(x)u'(x) dx \le 2 \| u \|_{L^2(a, \xi)} \| u' \|_{L^2(a, \xi)}
                           \le  \| u \|_{L^2(a, \xi)} ^2 +  \| u' \|_{L^2(a, \xi)}^2.
\]
In the same way, we have $\{ u(\xi) \}^2 \le  \| u \|_{L^2(\xi, b)} ^2 +  \| u' \|_{L^2(\xi, b)}^2$.
Combining the above two inequalities, (\eqho) holds.    \End

\propspace
{\bf \lemgood}. {\sl Let $\cY := C_{2\pi}^{\beta}(\R, L^2(\R) )$ and 
$\cZ := C_{2\pi}^{\beta}(\R, L^\infty(\R) )$.
Then, the following  hold.  
\n
(i) If $v \in \cY$ and  $w \in \cZ$  then $vw \in \cY$ and \,
$\|vw\|_\cY \le \|v\|_\cY \|w\|_\cZ$.
\n
(ii) If $v \in \cZ$ and  $w \in \cZ$  then $vw \in \cZ$ and \,
$\|vw\|_\cZ \le \|v\|_\cZ \|w\|_\cZ$.
}
\medn

	The proof  of {\lemgood}  is not difficult  and we leave it to the readers.

\propspace
{\bf \defper}. 
For a Banach space $E$ we define the Banach space
\[
C_{2\pi}(\R, E) := \{ u \in C(\R, E) \semicolon u(t + 2 \pi) = u(t) \moji{for} t \in \R \}
\]
with the norm $\displaystyle  \| u \|_{E, \,\text{per}} := \max_{t \in \R} \| u(t) \|_E$.  \End
     
\propspace
\n
{\bf \lemcont}.  {\sl Let $\cU = C^{\beta}_{2\pi}(\R, H^1(\R) )$ and  $\cZ$
be the same spaces as described in {\lemgood}. We assume  $f \in C^2(I, \R)$
and set $\displaystyle C_* :=  \max_{|r| \,\le\, d / \sqrt{2}} |f'(r)| +
                              \frac{d}{\sqrt{2}} \max_{|r| \,\le\, d/\sqrt{2}} |f''(r)|$.
Then  the following holds.
\medbreak
                 (i) If \,$u$, $v \in B_{\,\cU}(0 \semicolon d)$ then we have \, $f(u) - f(v) \in \cZ$ \, 
with the estimate 
\[
  \| f(u) - f(v) \|_\cZ \le C_* \| u - v \|_\cZ.      \tag{\eqau}
\]

               (ii)  If \,$u \in B_{\,\cU}(0 \semicolon d)$ then we have \, $f(u)  \in \cZ$ \, 
with the estimate: 
\[
           \| f(u)  \|_\cZ \le |f(0)| +  C_* \| u \|_\cZ.  
\]
}
 \medbreak
 
         {\it Proof}. (i) 
 Let $a$, $b \in (-d, d)$ and  $u$, $v \in B_{\,\cU}(0 \semicolon d)$.  We have
\[
f(a) - f(b) = \int_0^1  f' (\theta a + (1 - \theta)b) \, d \theta \cdot (a - b).    \tag{\eqax}
\]
It follows that
\[
f(u(x, t)) - f(v(x, t)) = g(x,t) \, \{ u(x, t) - v(x, t) \},     \tag{\eqg}
\]
where $\ds g(x,t) := \int_0^1 f' (\theta \, u(x, t) + (1 - \theta) v(x, t)) \, d \theta$.
By  {{\lemhone}, $g \in  C_{2\pi}(\R, L^\infty(\R) )$ and 
$\displaystyle  \| g(t)\|_{L^\infty(\R)} \le \max_{|r| \,\le\, d/\sqrt{2}} |f'(r)| \moji{for} t \in \R$.
It follows from (\eqax), {\lemhone} and Sobolev embedding  theorem that  for $x, s, t \in \R$ 
\begin{align*}
  & | g(x, t) - g(x,s) |   \\ 
= & \int^1_0 d \theta  \int^1_0 d \omega \,\,\,
     | f''(\omega \{ \theta u(x,t) + (1 - \theta)v(x,t) \} +(1 - \omega) \{ \theta u(x, s) + (1 - \theta)v(x, s) \}) | \\
   & \hspace{4.1 cm}    \cdot    \{ \theta || u(t) - u(s) ||_{L^\infty(\R)} +  
                                                                 (1 - \theta) || v(t) - v(s) ||_{L^\infty(\R)} \}    \\ 
\le \,&  \frac{1}{\sqrt{2}} \max_{|r| \,\le\, d/\sqrt{2}} |f''(r)| \, 
                 \int^1_0 d \theta \,\, \{ (\theta || u ||_\cU
                         +   (1 - \theta) || v ||_\cU  \} \,\, |t - s|^\beta  \\
\le \,& \frac{d}{\sqrt{2}} \, \max_{|r| \,\le\, d/\sqrt{2}} |f''(r)| \,\, |t - s|^\beta.
\end{align*}
It follows that
$ \displaystyle
\| g(t)  - g(s)\|_{L^\infty(\R)} \le  \frac{d}{\sqrt{2}} \max_{|r| \,\le\, d/\sqrt{2}} |f''(r)| \, |t - s|^\beta.
$
Therefore,  we have $g \in \cZ$ and 
$\ds  || g ||_\cZ \le C_*$.
By (\eqg) and {\lemgood} (ii), $f(u) - f(v) \in \cZ$ and (\eqau) holds.    
\medbreak
                    (ii) By (i) we have  $f(u) - f(0) \in \cZ$ and   $\| f(u) - f(0) \|_\cZ \le C_* \| u \|_\cZ$.
So, we obtain the desired result.     \End
 
\propspace
\n
{\bf \lemqua}. {\sl We assume (A). Then (\Hone-2) holds.}
\medn
{\it Proof.\,} 
Let $\cY$ and $\cZ$ be the same spaces as in \lemgood, and $\cU$ in \lemcont.
Let $(\lambda, u) \in \R \times B_X(0 \semicolon d)$.
Then we have $u_{xx} \in \cY$ and $u \in \cU$ with $u \in B_\cU (0 \semicolon d)$.
In view of {\lemcont} (ii), $\kappa_1(u) \in \cZ$. So, $\{ \kappa_1(u) - 1 \} u_{xx} \in \cY$ holds by 
{\lemgood} (i). In the similar way we verify that $\kappa_1'(u) u_x^2$, $u v^2$ and so on are 
elements of $\cY$. Therefore, (\HoneSecond) holds.
 \End

\propspace
         We define the map  
$\Psi : (\lambda, \ubold) \in \R \times B_X(0 \semicolon d) \longmapsto h(\lambda, \ubold) \in Y$.
\propspace
\n
{\bf \lemfs}. { \sl 
        We assume (A). Then (\Hone - 3) holds.
}
\medn
{\it Proof.\,} 
In view of {\lemqua}, 
we can redefine the maps $H$ and $h_0$ in (\eqhu)  \linebreak
as $H : B_X(0 \semicolon d) \to Y$,
$h_0 :  \R \times B_X(0 \semicolon d) \to Y$.

                  Let $\cX := C_{2\pi}^{1+\beta}(\R, L^2(\R) ) \cap C_{2\pi}^{\beta}(\R, H^2(\R) )$
 and $\cY := C_{2\pi}^{\beta}(\R, L^2(\R) )$.
Let $\kappa$ be the function satisfying the condition (A) with $\kappa_1$
replaced by $\kappa$.  We define the maps
$\varphi \colon u \in B_{\cX}(0 \semicolon d) \mapsto \kappa'(u)(u_x)^2 \in \cY$,
$\psi \colon u \in B_{\cX}(0 \semicolon d) \mapsto \{ \kappa(u) - 1\} u_{xx} \in \cY$,
$\gamma \colon u \in \cX \mapsto u^3 \in \cY$ and
$\omega \colon (u, v) \in X \mapsto u v^2 \in \cY$. 
In view of (\eqhu), it suffices to show that  the maps $\varphi, \psi, \gamma$ and $\omega$ 
are $C^2$ in order to prove (\Hone-3). We will show $\varphi \in C^2$ here   and leave
 the reader to prove $\psi, \gamma, \omega \in C^2$ since the proofs are similar.
         We verify 
\begin{gather*}
 D\varphi(u) v = \kappa''(u) (u_x)^2 v  + 2 \kappa'(u) u_x v_x        \tag{\eqdva},
\\
 D^2\varphi(u) v w = \kappa'''(u) (u_x)^2 v  w + 2 \kappa''(u) u_x (v_x w + v w_x) 
                                + 2 \kappa'(u)  v_x w_x      \tag{\eqdtv}                            
\end{gather*}
for  $u \in B_{\cX}(0 \semicolon d)$ and $v, w \in \cX$.
We denote by $\cL_2(E, F)$ the space of continuous bilinear maps from $E \times E \to F$
for the Banach spaces $E$ and $F$.
To prove the continuity of $D^2 \varphi$ it suffices to show
\[
\| D^2 \varphi(u_1) -  D^2 \varphi(u_2) \|_{\cL_2(\cX, \cY)} \le C \| u_1 - u_2 \|_{\cX}
                                                      \moji{for} u_1, u_2 \in B_{\cX}(0 \semicolon d),    \tag{\eqdt}
\]
where $C > 0$ is a constant independent of $u_1$ and $u_2$.
Actually, with respect to the first term in the right-hand side of (\eqdtv) it follows from 
Lemmas 5.1 -5.3 that 
for $u_1$, $u_2 \in B_{\cX}(0 \semicolon d)$ and $v, w \in \cX$  
\begin{align*}
& \| \kappa'''(u_1) (u_{1x})^2 v  w - \kappa'''(u_2) (u_{2x})^2 v  w \|_\cY     \\ 
\le\, &  \| \{ \kappa'''(u_1)  -  \kappa'''(u_2) \} (u_{1x})^2 v  w \| _\cY 
                                +  \| \kappa'''(u_2) \{ (u_{1x})^2 -  u_{2x})^2 \} v  w \|_\cY \\
\le\, & \|  \kappa'''(u_1)  -  \kappa'''(u_2)  \|_\cZ \, \| u_{1x} \|_\cZ^2  \,  \| v \|_\cZ  \, \| w \|_\cY \\
  &     \hspace{1cm}      +  \| \kappa'''(u_2) \|_\cZ  \, \| u_{1x}  + u_{2x} \|_\cZ \,
                                                                \| u_{1x} - u_{2x} \|_\cZ  \,\| v \|_\cZ \, \| w \|_\cY  \\
\le\, & \left(  \frac{C_1}{4} \,  \| u_1  \|_{\cX}^2 +   \frac{C_2}{2\sqrt{2}} \, \| u_1 + u_2  \|_{\cX} \right)
                       \| u_1 - u_2 \|_{\cX} \| v \|_{\cX} \| w \|_{\cX}  \\
 \le \,& \,  C\, \| u_1 - u_2 \|_{\cX} \| v \|_{\cX} \| w \|_{\cX},                                                  
\end{align*}
where we set $\displaystyle C_1\mycolon =  \max_{|r| \,\le\, d / \sqrt{2}} |\kappa^{(4)} (r)| +
                                 \dfrac{d}{\sqrt{2}}  \max_{|r| \,\le\, d / \sqrt{2}} |\kappa^{(5)} (r)|$,
$\displaystyle C_2\mycolon = |\kappa'''(0)| + \dfrac{C_1d}{\sqrt{2}}$
\linebreak 
  and $C \mycolon = \dfrac{C_1 d^2}{4} +  \dfrac{C_2 d}{\sqrt{2}}$.
In the similar way we can estimate the other terms in the right-hand side of (\eqdtv).
So, (\eqdt) holds.       \End

\propspace
{\bf \lemaa}. { \sl 
We assume (A). Then (\Hone - 4) holds.
}
\medn
{\it Proof.\,} Clearly, $h(\lambda, {\bm 0})=  {\bm 0}$ for $\lambda \in \R$.
Let $\ubold \in B_U(0 \semicolon d)$ and $\hbold = (h, k) \in U$. 
Then, we verify 
\begin{gather*}
D_{\bm u} h_0 (0, {\bm u}) {\bm h}  = 
\begin{pmatrix}
- 3 u^2h - v^2h + 2uvk    \\ 
-2 uvh - u^2k - 3v^2k
\end{pmatrix},     \tag{\eqdh}   \\
DH(\ubold)\hbold =
\begin{pmatrix}
\{ \kappa_1' (u) u_x h +  (\kappa_1(u) - 1) h_x \} _x      \\ 
\{\kappa_2' (v) v_x k +  (\kappa_2(v) - 1) k_x \} _x
\end{pmatrix}.                                               
\end{gather*}
It follows that
$D_{\bm u} h_0 (0, {\bm 0}) = 0$  and $DH({\bm 0})  =  0$.   
So,   $D_{\bm u} h(0, {\bm 0})= DH({\bm 0}) +  D_{\bm u} h_0(0, {\bm 0}) =0$.    \End

\newchapter

\centerline{\bf References}
\medn
[\refAma] H. Amann, Hopf bifurcation in quasilinear reaction--diffusion systems,
Delay \linebreak
Differential Equations and Dynamical Systems, Lecture Notes in Mathematics 1475 (1991) 53--63.
\medn
[\refABB] W. Arendt, C.Batty and S. Bu, Fourier
multipliers for H{\"o}lder continuous functions and maximal regularity,
Studia mathematica 160 (1) 23-51, 2004.
\medn
[\refBKST] T. Brand, M. Kunze, G. Schneider, T. Seelbach, Hopf bifurcation and 
exchange of stability in diffusive media, Arch. Rat. Mech. Anal. 171 (2004) 263--296.
\medn
[\refCRbif] M.G. Crandall and P.H. Rabinowitz,  Bifurcation from  
simple eigenvalues, J. Func. Anal. 8 (1971) 321--340. 
\medn
[\refCRhopf] M.G. Crandall and P.H. Rabinowitz, 
The Hopf bifurcation theorem in
 infinite dimensions,
 Arch. Rat. Mech. Anal. {67} (1977) 53--72.
 \medn
[\refGMW] D. Gomez, L. Mei, J. Wei,
Stable and unstable periodic spiky solutions for the Gray-Scott system and
 the Schnakenberg system, J. Dynam. Differential Equations 32 (2020), 441--481. 
 \medn
[\refKsym] T. Kawanago, A symmetry-breaking bifurcation 
theorem and some related theorems
 applicable to maps having unbounded derivatives,
Japan J. Indust. Appl. Math. 21 (2004) 57--74. 
Corrigendum to this paper:  
Japan J. Indust. Appl. Math. 22 (2005) 147. 
\medn
[\refKcom] T. Kawanago,  Computer assisted proof to symmetry-breaking bifurcation phenomena in nonlinear vibration,
Japan J. Indust. Appl. Math. 21 (2004) 75--108. 
 \medn
[\refKsome] T. Kawanago, Codimension-$m$ bifurcation theorems 
applicable to the numerical\linebreak 
verification methods, Advances in Num. Analysis, 
vol. 2013 (2013), Article ID 420897. 
\medn
[\refKhi] T. Kawanago, The Hopf bifurcation theorem in Hilbert spaces for abstract semilinear 
equations, J. Dynamics and Differential Equations, 35 (2023) 2677–2690.
 \medn
[\refKiel]  H. Kielh{\"o}fer, Bifurcation theory. An introduction with 
applications to partial \linebreak 
differential equations. Second edition. 
Applied Mathematical Sciences, 156. Springer, New York, 2012. 
\medn
[\refLiZY] H. Li,  X. Zhao and W. Yan, Bifurcation of time-periodic solutions for the 
incompressible flow of nematic liquid crystals in three dimension, Adv. Nonlinear Anal. 9 (2020) 1315--1332.
 \medn
[\refliu] Z. Liu, P. Magal and S. Ruan, Hopf bifurcation for 
non-densely defined Cauchy problems,
Z. Angew. Math. Phys. 62 (2011) 191--222.
\medn
[\refMS] A. Melcher and G. Schneider,
A Hopf-bifurcation theorem for the vorticity formulation of the Navier-Stokes equations in
$\R^3$,  Comm. Partial Differential Equations 33 (2008), no. 4-6, 772--783.

\end{document}